\documentstyle[12pt,epsf]{article}

\newcommand{\rr}{{\rm{I \! R}}}
\newcommand{\ep}{\varepsilon}
\newcommand{\gam}{\gamma}
\newcommand{\suml}{\sum\limits}
\newcommand{\al}{\alpha}
\newcommand{\ov}{\overline}
\newcommand{\widt}{\widetilde}
\newcommand{\ds}{\displaystyle}
\newcommand{\f}{\frac}
\newcommand{\diag}{{\rm{diag}\,}}
\newcommand{\fii}{\varphi}
\newcommand{\lam}{\lambda}

\newcommand{\ff}{{\rm{I \! F}}}
\newcommand{\ver}{{\rm {ver}}\,}
\newcommand{\Om}{\Omega}
\newcommand{\De}{\Delta}
\newcommand{\om}{\omega}
\newcommand{\de}{\delta}
\newcommand{\p}{\partial}
\newcommand{\calO}{O^{\!\!\!\!\smile}}
\newcommand{\grad}{{\rm{grad}\,}}
\newcommand{\Rem}{{\rm{Re}\,}}
\newcommand{\ca}{{\cal A}}
\newcommand{\th}{\theta}
\newcommand{\cb}{{\cal B}}
\newcommand{\cv}{{\cal V}}
\newcommand{\loc}{{\rm loc}}
\newcommand{\cc}{{\rm{{\footnotesize{l}}\!\!\! C}}}

\title{\bf DISSIPATIVE MECHANICAL SYSTEMS WITH DELAY}
\author{Ion D. ALBU, Mihaela NEAM\c TU, Dumitru OPRI\c S}
\date{}

\begin{document}
\maketitle

\begin{quote}{\small{{\bf Abstract.}
The idea of dissipative mechanical system with delay is proposed. The paper studies
the phenomenon of dissipation with delay for Euler-Poincar\'e systems on Lie
algebras or equivalently, for Lie-Poisson systems on the duals of Lie algebras. The
study was suggested by the work [2] and it is ended with a discussion regarding the
stability and the Hopf bifurcations for the free rigid body with delay.}}
\end{quote}

\noindent{\small{{\it Keywords:} delay differential equation,
dissipation with delay, stability, Hopf bifurcation.

\noindent{\it 2000 AMS Mathematics Subject Classification:} 34K13, 34K20, 37G15,
37J15, 53D20.}}

\section*{\normalsize\bf Introduction}

\hspace{0.6cm} In many applications one assumes the system under consideration is governed by a principle of
causality; that is, the future state of the system is independent of the past states and is determined solely by
the present. If it is also assumed that the system is governed by an equation involving the state and the rate of
change of the state, then, generally, one considers either ordinary or partial differential equations. However,
under a closer scrutiny, it becomes apparent that the principle of causality is often only a first approximation
to the true situation and that a more realistic model would include some of the past states of the system. The
simplest type of past dependence is through the state variable and not the derivative of the state variable, the
so--called retarded functional differential equations or retarded difference equations or systems with delay.
Systems with delay are studied in many biological research topics, as well as in several branches of engineering,
in neural networks, in economics, optimal production decision for an oligopoly with information lags ([9], [11],
[12] etc.).

Functional differential equations on finite--dimensional manifolds are considered in [5], where the topological
properties of the global attractor of an retarded functional differential equation in terms of limit capacity and
Hausdorff dimension are presented.

The purpose of this paper is to study the phenomenon of dissipation with delay for Euler--Poincar\' e systems on
Lie algebras or equivalently, for Lie--Poisson systems on the duals of Lie algebras. The dissipation without delay
inducing instabilities for Euler--Poincar\'e systems is studied in [2]. The dissipation with delay that we
construct has the essential feature: the energy is dissipated, but the angular momentum is not. In the context of
Euler--Poincar\'e or Lie--Poisson systems this means that the coadjoint orbits remain invariant but on them the
energy is decreasing along orbits. It is interesting the geometry behind the construction of the nonlinear
dissipative term with delay, which has a Brockett with delay double bracket form. In fact, this form is well
adapted to the study of dissipation with delay on Lie groups since it was constructed as a gradient system and it
is well known in other contexts that this formalism plays an important role in the study of integrable systems.

The general equations of motion for dissipative systems with delay
that we consider have the following form:
$$\dot F=\{F,H\}-\left\{\left\{\widt F,\widt H\right\}\right\}_d$$
where $H$ is the total energy of the system, $\{F,H\}$ is a skew
symmetric bracket which is a Poisson bracket in the usual sense
and where $\left\{\left\{\widt F,\widt H\right\}\right\}_d$ is a
{\it symmetric bracket with delay}.

The outline of this paper is as follows. In section 1 the
functional diffe\-ren\-tial equations on manifolds are introduced.
In section 2 a concrete idea of dissipative mechanism with delay
is presented. The same formalism can be applied to other systems
as well. Some of the basic and essential facts about dissipative
mechanical systems with delay are described in section 3. In
section 4 we study Lagrangian systems that are invariant under a
group action and we shall add to them, in the sense of section 3,
dissipative fields with delay that are equivariant. It is given by
a necessary and sufficient condition that the integral curves of
the vector field $Z+Y$, for a vertical $G$--invariant vector field
$Y$ on $TQ\times TQ$ and a Lagrangian vector field $Z$ of a
$G$--invariant Lagrange function, preserve the inverse images of
the coadjoint orbits in ${\bf g}^*$ by the momentum map $J$. In
section 5 the Euler--Poincar\'e and Lie--Poisson equations with
delay are given and the double bracket with delay is defined; in
the case ``without delay" this is defined in [2]. In section 6 is
presented the free rigid body with delay if the components
$(I_1,I_2,I_3)$ of the moment of the inertial tensor satisfy the
conditions $I_1>I_2$, $I_1>I_3$. For the equilibrium state
$\Om_1=(m/I_1,0,0)^T$, $m\not=0$, the value $\tau_c$ of the delay
for which $\Om_1$ is asymptoticaly stable and the value $\tau_0$
of the delay for which there is a Hopf bifurcation are determined.
They are also determined the local center manifold and the
quantities $C_1(0)$, $\mu_2, T_2,\beta_2$ sketching out the
direction of the Hopf bifurcation, the stability and the period of
the bifurcating periodic solutions. These quantities are
calculated for fixed values of $\al$ and $m$. Finally, some
conclusions are drawn and further research directions are
discussed in the last section.

\section*{\normalsize\bf{1. Functional differential equations on
manifolds}}

\hspace{0.6cm} In this section, we begin with examples that will
serve as a motivation for the consideration of functional
differential equations on manifolds.

{\bf Example 1.} For any constant $c$, the scalar equation
$$\dot q (t)=c\sin (q(t-1))\eqno(1.1)$$
can be considered as an retarded functional differential equations
(RFDE) on the circle
$S^1=\left\{x\in\rr^2,~(x^1)^2+(x^2)^2=1\right\}$ by considering
$q$ as an angle variable only determined up to a multiple of
$2\pi$.

{\bf Example 2.} If $b,c$ are constants then we write the
second--order RFDE
$$\ddot q (t)+b\dot q(t)=c\sin (q(t-1))\eqno(1.2)$$
as a system of first--order RFDE
$$\begin{array}{l}
\vspace{0.1cm}
\dot q{}^1(t)=q^2(t),\\
\dot q{}^2(t)=c\sin (q^1(t-1))-bq^2(t).\end{array}\eqno(1.3)$$ By
considering $q^1$ as an angle variable only determined up to a
multiple $2\pi$, the equation (1.3) is a RFDE on the cylinder
$S^1\times\rr$. We remark that we can take the space of initial
data for the solution $(q^1,q^2)$ of (1.3) as
$C([-1,0],S^1)\times\rr$.

The simplest type of past dependence in a differential equation is
that in which the past dependence is through the state variable
and not the derivative of the state variable, the so - called
retarded functional differential equations or retarded
differential difference equations. For a discussion of the
physical applications of the differential difference equation
$$\dot q(t)=F(t,q(t), q(t-\tau),$$
to control problems, see [5], [6].

{\bf Example 3.} The equations involved in the study of vibrating masses attached to an elastic bar are
$$\begin{array}{l}
\vspace{0.1cm}
\ddot q{}^1(t)+\omega_1^2q^1(t)=\ep f_1(q^1(t),\dot q{}^1(t), q^2(t),\dot q{}^2(t))
+\gam_1\ddot q{}^2(t-\tau),\\
\ddot q{}^2(t)+\omega^2_2 q^2(t)=\ep f_2\left(q^1(t), \dot
q{}^1(t), q^2(t),\dot q{}^2(t)\right)+\gam_2\ddot
q{}^1(t-\tau).\end{array}\eqno(1.4)$$

{\bf Example 4.} Let
$S^2=\left\{(q^1,q^2,q^2)\in\rr^3,~(q^1)^2+(q^2)^2+(q^3)^2=1\right\}$
and consider the following system of RFDE:
$$\begin{array}{l}
\vspace{0.1cm}
\dot q{}^1(t)=-q^1(t-\tau)q^2(t)-q^3(t),\\
\vspace{0.1cm}
\dot q{}^2(t)=q^1(t-\tau)q^1(t)-q^3(t),\\
\dot q{}^3(t)=q^1(t)+q^2(t).\end{array}\eqno(1.5)$$ If
$(q^1(t),q^2(t),q^3(t))$ is a solution of the equation (1.5), it
is easy to see that
$$\sum_{i=1}^3 q^i(t)\dot q{}^i(t)=0$$
for all $t\ge 0$, $\forall\tau$. As a consequence, for $t\ge 0$,
$\suml_{i=1}^3\left(q^i(t)\right)^2=a^2$, $a$ constant. Thus, if
an initial condition $\phi=(\phi^1,\phi^2,\phi^3)$ satisfies
$\phi(\theta)\in S^2$ for all $\theta\in[-\tau,0]$, we conclude
that the solution $(q^1,q^2,q^3)(t;\phi)\in S^2$ for all $t\ge 0$.
With this remark, we can define an RFDE on $S^2$ by the map
$$F:\phi=(\phi^1,\phi^2,\phi^3)\in C([-\tau,0],S^2)\to F(\phi)\in TS^2$$
where $F(\phi)$ is the tangent vector to $S^2$ at the point
$\phi(0)$ given by
$$F(\phi)=(-\phi^1(-\tau)\phi^2(0)-\phi^3(0), \phi^1(-\tau)\phi^1(0)-\phi^3(0),\phi^1(0)+\phi^2(0)).$$
We now formalize the notions in these examples to obtain a RFDE on
a $n$--dimensional manifold. Roughly speaking, a RFDE on a
manifold $Q$ is a function $F$ mapping each continuous path $\phi$
lying on $Q$, $\phi\in C([-\tau,0], Q]$, into a vector tangent
$F(\phi)$ to $Q$ at the point $\phi(0)\in Q$.

Let $Q$ be a separable $C^\infty$ finite $n$--dimensional manifold
(configuration manifold), $I=[-\tau,0]$, $\tau\ge 0$, and $C(I,Q)$
the totality of the continuous maps $\phi$ of $I$ into $Q$. The
space $C(I,Q)$ is separable and is a $C^\infty$--manifold modeled
on a separable Banach space.

If $\rho:C(I,Q)\to Q$ is the evaluation map, $\rho(\phi)=\phi(0)$,
then $\rho$ is $C^\infty$, and for each $q\in Q$, $\rho^{-1}(q)$
is a closed submanifold of $C(I,Q)$ of codimension $n=\dim Q$.

A {\it retarded functional differential equation} (RFDE) on $Q$ is
a continuous function $F:C(I,Q)\to TQ$ such that $\pi_{TQ}\circ
F=\rho$. If we want to emphasize the function $F$ defining the
RFDE, we write RFDE$(F)$.

A {\it solution} of RFDE$(F)$ is defined in the obvious way, namely, as a con\-ti\-nuous function
$q:[-\tau,\al)\to Q$, $\al>0$, such that $\dot q(t)$ exists and is continuous for $t\in [0,\al)$ and $(q(t),\dot
q(t))=F(q_t)$, for $t\in [0,\al)$ where $q_t(\theta)=q(t+\theta)$, $\theta\in [-\tau,0]$. Locally, if
$F(\phi)=(\phi(0), f(\phi))$ for an appropiate function $f$, then this is equivalent to $\dot q(t)=f(q_t)$.

The basis theory of existence, uniqueness and continuous dependence on initial data for general RFDE on manifold
is the same as the theory when $Q=\rr^n$ [4].

Any $C^k$--vector field on $Q$ defines a $C^k$--RFDE on $Q$. In
fact, if $X:Q\to TQ$ is a $C^k$--vector field on $Q$, it is easy
to see that $F=X\circ\rho$ is a $C^k$--RFDE on $Q$.

To show that the equation considered in Example 2 is an RFDE
according to the definition, we need the concept of a second order
RFDE on $Q$. Let $\ov F: C(I, TQ)\to TQ\times TQ$ be a continuous
function that locally has the representation
$$\ov F(\phi,\psi)=\left(\left(\phi(0),\psi(0)\right), \left(\psi(0), f(\phi,\psi)\right)\right).$$
The solution $(x(t), y(t))$ of the RFDE$(\ov F)$ on $TQ$ satisfies
the equations
$$\dot q(t)=y(t),\quad \dot y(t)=f(x_t, y_t)$$
where $q(t)\in Q$. If it is possible to perform the
differentiations, then we obtain the second--order equation
$$\ddot q(t)=f(q_t,\dot q_t).\eqno(1.6)$$

If we now return to Example 2, we see that the formulation
requires that we consider initial data in the space
$C(I,S^1)\times C(I,\rr)$. However, this does not affect the
dynamics since the solution will be in space $C(I,S^1)\times\rr$
after one unit of time.

Let $g:Q\times Q\to TQ$ be such that $(\pi_{TQ}\circ g)(q,\widt
q)=q$ and let $d:C(I,Q)\to Q\times Q$ be defined by
$d(\phi)=(\phi(0),\phi(-\tau))$. The function $F=g\circ d$ is a
RFDE on $Q$ which can be written locally as
$$\dot q(t)=\ov g(q(t), q(t-\tau))\eqno(1.7)$$
where $g(q,\widt q)=(q,\ov g(q,\widt q))$. The equations (1.7) is
a {\it delay differential equation} (DDE) on $Q$.

A $C^k$--{\it vector field with delay} is given by the mappings
$g:Q\times Q\to T^*Q$, such that $(\pi_{T^*Q}g)(q_1,q_2)=q_1$, and
$d:C(I,Q)\to Q\times Q$, defined by
$d(\phi)=(\phi(0),\phi(-\tau))$. The function $G=g\circ d$ is a
covector field with delay on $Q$ that can be written locally as
$$\omega(t)=f_i(q(t), q(t-\tau))dq^i\eqno(1.8)$$
where $g(q,\widt q)=\left(q, f(q,\widt q)\right)$.

The topological properties of the RFDE on manifolds are discussed in [5].

\section*{{\normalsize\bf 2. Motivating examples}}

\hspace{0.6cm} To get a concrete idea of the type of dissipative
mechanisms with delay we have in mind, we now give a simple
example of it for perhaps the most basic of Euler--Poincar\'e or
Lie--Poisson systems, namely the rigid body. Here, the Lie algebra
in question is that of the rotation group; that is, the Euclidean
space $\rr^3$ interpreted as the space of body angular velocities
$\Omega$ equipped with the cross product as the Lie bracket. On
this space, we put the standard kinetic energy Lagrangian
$L(\Omega)=\ds\f{1}{2}(I\cdot\Omega)$ [where $I=\diag
(I_1,I_2,I_3)$ is the moment of inertial tensor] so that the
general Euler--Poincar\'e equations become the standard rigid body
equations for a freely spinning rigid body:
$$I\dot\Omega=(I\Omega)\times\Omega\eqno(2.1)$$
or, in terms of the body angular momentum $M=I\Omega$,
$$\dot M=M\times\Omega.\eqno(2.2)$$
In this case, the energy equals the Lagrangian;
$E(\Omega)=L(\Omega)$ and the energy is conserved by the solutions
of (2.1). Now we modify the equations by adding a term with delay
$$\dot M=M\times\Omega+\al M\times \left(\widt M\times \Omega\right)\eqno(2.3)$$
where $\al$ is a positive constant, $\widt M(t)=M(t-\tau)$, for
$t\ge 0$, $\tau\ge 0$ and $\widt M(t)=\fii(t)$, for $t\in
[-\tau,0]$.

A related example is the Laudau--Lifschitz equations with delay
for the magnetization vector $M$ in a given magnetic field $B$
$$\dot M=\gam M\times B+\f{\lam}{\|M\|^2\cdot\cos\theta}\left(
M\times\left(\widt M\times\widt B\right)\right)\eqno(2.4)$$ where
$\gam$ is the magneto--mechanical ratio, $\lam$ is the damping
coefficient due to domain walls and $\theta$ is the angle between
$M$ and $\widt M$.

One checks in each case that the addition of the dissipative term
with delay has a member of interesting properties. First of all,
this dissipation with delay is derivable from a $SO(3)$--
invariant force field. However, it is induced by a dissipation
function with delay in the following restricted sense. It is a
gradient when is restricted to each momentum sphere (coadjoint
orbit) where each sphere carries a special metric (later to be
called the normal metric). Namely, the extra dissipative term with
delay in (2.3) equals the negative gradient of the Hamiltonian
with respect to the following metric on the sphere. Each vector
$v\in\rr^3$ can be orthogonally decomposed with respect to the
standard metric on $\rr^3$ into a component tangent at $M$ to the
sphere $\|M\|^2=c^2$ and a component on $\widt M$, where
$\left\|\widt M\right\|=c^2$, $\widt M\not= M$:
$$v=\f{M\cdot v}{c^2\cos\theta}\widt
M-\f{1}{c^2\cos\theta}\left[M\times\left (\widt M\times
v\right)\right]\eqno(2.5)$$ where $\theta$ is the angle from $M$
and $\widt M$. The metric on the sphere is chose to be
$(c^2\cos\theta)^{-2}\al$ times the standard inner product of the
components tangent to the sphere in the case of the rigid body
model with delay and just $\lam$ times the standard metric in the
case of the Laudau--Lifschitz equations with delay.

Secondly, the dissipation with delay to the equations has the
obvious form of a repeated Lie bracket, i.e. a double bracket, and
it has the properties that the conservation law
$$\f{d}{dt}\|M\|^2=0\eqno(2.6)$$
is preserved by the dissipation with delay (since the extra force
is orthogonal to $M$) and the energy is strictly monotone except
at relative equilibria. In fact, we have:
$$\f{d}{dt}E=-\al\left\|\widt M\times\widt\Omega\right\|^2\eqno(2.7)$$
for the rigid body and
$$\f{d}{dt}E=-\f{\lam}{\|M\|^2\cos\theta}\left\|\widt M\times\widt B\right\|^2\eqno(2.8)$$
in the case of the Laudau--Lifschitz equations, so that the
trajectories on the angular momentum sphere converge to the
minimum (for $\al$ and $\lam$ positive) of the energy restricted
to the sphere, apart from the set of measure zero consisting of
orbits that are relative equilibria or are the stable manifolds of
the perturbed saddle point.

Another interesting feature of the dissipations with delay is that
they can be derived from a bracket in the same way that the
Hamiltonian equations can be derived from a skew symmetric Poisson
bracket. For the case of the rigid body with delay, this bracket
is
$$\left\{\{F,K\}\right\}=\al(M\times\nabla\widt F)\cdot \left(\widt M\times\nabla\widt K\right).$$

As we have already indicated, the same formalism can be applied to
other systems as well. In fact, later in the paper we develop an
abstract constructions for dissipative with delay terms with the
same general properties as the above examples.

\section*{\normalsize\bf 3. Dissipative systems with delay}

\hspace{0.6cm} For later use, it will be useful some of the basic
and essential facts about dissipative mechanical systems with
delay. Let $Q$ be a manifold, $L:TQ\to\rr$ be a smooth function
and let $\pi:TQ\to Q$ be the tangent bundle projection. Let $\ff
L:TQ\to T^*Q$ be the {\it fibre derivative} of $L$; recall that it
is defined by
$$<\ff L(v), w>=\f{d}{d\ep}\Bigl|_{\ep=0} L(v+\ep w).\eqno(3.1)$$
Where $<,>$ denotes the pairing between the tangent and cotangent
spaces. We also recall that the {\it vertical lift} of a vector
$w\in T_q Q$ along $v\in T_qQ$ is defined by:
$$\ver_v(w)=\f{d}{d\ep}\Bigl|_{\ep=0}(v+\ep w)\in T_v(TQ).\eqno(3.2)$$
The action and energy of $L$ are defined by
$$A(v)=<\ff L(v),v>\eqno(3.3)$$
and
$$E(v)=A(v)-L(v).\eqno(3.4)$$

Let $\Om_L=(\ff L)^*\Om$ denote the pull back of the canonical
sympletic form on $T^*Q$ by the fibre derivative of $L$.

A vector field $Z$ on $TQ$ is called a {\it Lagrangian vector
field} of $L$ if
$$i_Z\Om_L=dE.\eqno(3.5)$$
In this generality, $Z$ need not exist, nor be unique. However, we shall assume throughout that $Z$ is a second
order equation; that is $T{\pi\circ Z}$ is the identity on $TQ$. A second order equation is a Lagrangian vector
field if and only if the Euler--Lagrange equations hold in local charts. We note that, by skew symmetry of
$\Om_L$, the energy is always conserved; that is, $E$ is constant along an integral curve of $Z$. We also recall
that the Lagrangian is called {\it regular} if $\Om_L$ is a (weak) sympletic form; it is nondegenerate. This is
equivalent to the second fibre derivative of the Lagrangian being, in local charts, also weakly nondegenerate. In
the regular case, if the Lagrangian vector field exists, it is unique, and is given by the Hamiltonian vector
field with energy $E$ relative to the sympletic form $\Om_L$. If, in addition, the fibre derivative is a global
diffeomorphism, then $Z$ is the pull back by the fibre derivative of the Hamiltonian vector field on the cotangent
bundle with Hamiltonian  $H=E\circ (\ff L)^{-1}$. It is well known how one can pass back and forth between the
Hamiltonian and Lagrangian pictures in the hyperregular case [2].

Consider a general Lagrangian vector field $Z$ for a (not
necessarily regular) Lagrangian on $TQ$. A map $Y:TQ\times TQ\to
TQ$ is called a {\it dissipative vector field with delay} if is
vertical, i.e. $T\pi\circ Y=0$ and if at each point of $TQ\times
TQ$
$$<dE, Y>~\le~ 0.\eqno(3.6)$$ If the inequality is
pointwise strict at each nonzero $v\in TQ$, $\widt v\in TQ$,
$v\not=\widt v$, then we say that the  map $Y$ is dissipative. A
{\it dissipative Lagrangian system} on $TQ$ is a vector field of
the form $X=Z+Y$, where $Z$ is a (second order) Lagrangian vector
field and $Y$ is a dissipative vector field with delay. Define the
1--form $\De^Y$ on $TQ\times TQ$ by
$$\De^Y=-i_Y\Om_L\eqno(3.7)$$
and the force field with delay $F^Y:TQ\times TQ\to T^*Q$ given by
$$<F^Y\left(v,\widt v\right), w>=\De^Y(v,v)\cdot V_v=-\Om_L(v)\left(Y\left(v,\widt v\right), V_v\right)\eqno(3.8)$$
where $T\pi(V_v)=w$, and $V_v\in T_v(TQ)$.

{\bf Proposition 3.1.} {\it A vertical vector field $Y:TQ\times
TQ\to TQ$ is dissipative with delay if and only if the induced
force field with delay $F^Y$ satisfies $<F^Y\left(v,\widt
v\right), v>\,<0$ for all nonzero $v\in TQ$ ($\le 0$ for the
weakly dissipative with delay).}

{\bf Proof.} Let $Y$ be a vertical vector field $Y:TQ\times TQ\to
TQ$ with $T\pi\circ Y=0$, $\De^Y$ the form on $TQ\times TQ$ given
by (3.7) and $F^Y$ the force field with delay given by (3.8). If
$Z$ denotes the Lagrangian system defined by $L$, we get
$$\begin{array}{lll}
\vspace{0.1cm}
(dE\cdot Y)\left(v,\widt v\right)&=&(i_Z\Om_L)(Y)\left(v,\widt v\right)=\Om_L(Z,Y)\left(v,\widt v\right)=\\
\vspace{0.1cm}&=& -\Om_L(v)\left(Y\left(v,\widt v\right), Z(v)\right)=\\
&=&< F^Y\left(v,\widt v), T_v\pi(Z(v)\right)>=<F^Y\left(v,\widt
v\right),v>,\end{array}$$ since $Z$ is a second--order equation.
We conclude that $<dE,Y>\,<0$ if and only if $<F^Y\left(v,\widt
v\right),v>\,<0$, for all $\left(v,\widt v\right)\in TQ\times TQ$,
$v\not=\widt v$, which gives the result.

Treating $\De^Y$ as the exterior force with delay $1$--form acting
on a mechanical system with a Lagrangian $L$, we now shall write
the governing equation of motion. The basic principle is of course
the Lagrange--d'{}Alembert principle.

The {\it Lagrangian force} associated with a given Lagrangian $L$
and a given second--order vector field $X$ is the horizontal
$1$--form on $TQ$ defined by
$$\phi_L(X)=i_X\Om_L-dE.\eqno(3.9)$$
Given a horizontal $1$--form $\om$ (referred to $1$--form as the
exterior force with delay), the local Lagrange--d'{}Alembert
principle states that
$$\phi_L(X)+\om=0.\eqno(3.10)$$

It is easy to check that $\phi_L(X)$ is indeed horizontal if $X$
is of second order. Conversely, if $L$ is regular and if
$\phi_L(X)$ is horizontal, then $X$ is of second order. One can
also formulate an equivalent principle in variational form.

Given a Lagrangian $L$ and a force field with delay (as defined in
Proposition 1) the {\it integral Lagrange--d'{}Alembert principle
with delay} for a curve $q(t)$ in $Q$ is
$$\de\int_a^b L\left(q(t),\dot q(t)\right)dt+\int_a^b F\left(\left(q(t),\dot q(t)\right),
 \left(\widt q(t),\dot{\widt q}(t)\right)\right)\cdot \de q dt=0,\eqno(3.11)$$
where the variation is given by the usual expression
$$\de\int_a^b L\left(q(t),\dot q(t)\right)dt=\int_a^b\left(\f{\p L}{\p q^i}-\f{d}{dt}\left(\f{\p L}{\p \dot q{}^i}\right)\right)\de q^i dt\eqno(3.12)$$
for a given variation $\de q$ (vanishing at the endpoints).

The two forms of the Lagrange--d'{}Alembert principle are
equivalent. This follows from the fact that both give the
Euler--Lagrange equations with for\-cing with delay in local
coordinates (provided that $Z$ is of second order). We shall see
this in the following development.

{\bf Proposition 3.2.} {\it Let $\om$ be the delay exterior force
$1$--form associated to a vertical vector field with delay $Y$,
i.e. $\om=\De^Y=-i_Y\Om_L$. Then $X=Z+Y$ satisfies the local
Lagrange-d'{}Alembert principle with delay. Conversely, if, in
addition, $L$ is regular, the only second-order vector field $X$
satisfying the local Lagrange--d'{}Alembert principle with delay
is $X=Z+Y$.}

{\bf Proof.} For the first part, the equality $\phi_L(X)+\om=0$ is
a simple ve\-ri\-fi\-cation. For the converse, we already know
that $X$ is a solution and the uniqueness is guaranteed by
regularity.

To develop the differential equations associated to $X=Z+Y$, we
take $\om=-i_Y \Om_L$ and note that, in a coordinate chart,
$Y\left(\left(q, \dot q\right)\left(\widt q,\dot{\widt
q}\right)\right)=Y^i\left(\left(q,\dot q\right),\left(\widt
q,\dot{\widt q}\right)\right)\ds\f{\p}{\p q^i}$ and the equation
(3.10) is given by
$$\f{d}{dt}\left(\f{\p L}{\p\dot q{}^i}\right)-\f{\p L}{\p q^i}=\f{\p^2 L}{\p\dot q{}^i\p\dot q{}^j}Y^j\left((q^k,\dot q{}^k),\left(\widt q{}^k,\dot{\widt q}{}^k\right)\right).\eqno(3.13)$$
The force $1$--form with delay $\De^Y$ is therefore given by
$$\De^Y\left(q^k,\dot q{}^k),\left(\widt q{}^k,\dot{\widt q}{}^k\right)\right)=\f{\p^2L}{\p\dot q{}^i\p\dot q{}^j} Y^j\left((q^k, \dot q{}^k),\left(\widt q{}^k,\dot{\widt q}{}^k\right)\right)dq^i\eqno(3.14)$$
and the corresponding force field with delay is given by
$$F_i^Y=\left(q^k,\f{\p^2L}{\p\dot q{}^i\p\dot q{}^j}Y^j\left((q^k,\dot q{}^k),
\left(\widt q{}^k,\dot{\widt
q}{}^k\right)\right)\right).\eqno(3.15)$$

Thus, the condition for an integral curve takes the form of the
Euler--Lagrange equations with force with delay
$$
\f{d}{dt}\left(\f{\p L}{\p\dot q{}^i}\right)-\f{\p L}{\p
q^i}=F_i^Y\left((q^k,\dot q{}^k),\left(\widt q{}^k,\dot{\widt
q}{}^k\right)\right).\eqno(3.16)$$

Since the integral Lagrange--d'{}Alembert principle with delay
gives the same equations, it follows that the two principles are
equivalent. From now on, we shall refer to either one as simple
the Lagrange--d'{}Almbert principle with delay.

{\bf Example 1.} The inertial neuron with time delay is described
by the equations
$$\ddot q{}^i=-a\dot q{}^i-bq^i+cf(q^i-h\widt q{}^i)+d\sum_{j=1\atop{j\not=i}}^n f
\left(q^j-h\widt q{}^j\right),\quad i=\ov{1,n},\eqno(3.17)$$ where $a,b,c,d>0$, $h\ge 0$ are constants and $\widt
q{}^i(t)=q^i(t-\tau)$ is the time delay. For $n=1$ the model is discussed in [1]. For $L\left(\dot q,
q\right)=\ds\f{1}{2}\ds\sum_{i=1}^n \left(\dot q{}^i\right)^2-b\ds\sum_{i=1}^n q^i$ and the force field with delay
given by
$$Y^i\left(\dot q,\widt q)=-a\dot q{}^i+cf(q^i-h\widt
q{}^i\right)+d\sum_{j=1\atop{j\not=i}}^n
 f\left(q^j-h\widt q{}^j\right),\eqno(3.18)$$
the equations (3.17) are the Euler--Lagrange equations with force
with delay. The model for $n=2$ is analyzed in [9].   The force
field with delay (3.18) is not dissipative.

{\bf Example 2.} The simplest mechanical model of the regenerative
machine tool vibration in the case of the so--called orthogonal
cutting is given by the equation
$$\ddot q+2k\al\dot q+\al^2 q=\f{1}{m}f\left(\dot q,\beta\right),\eqno(3.19)$$
where $\al=\sqrt{s/m}$ is the natural angular frequency of the undamped free oscillating system, $k$ is the
so--called relative damping factor and $f\left(\dot q,\beta\right)$ is the cutting force as a function of
technological parameters and as a function of the chip thickness $\beta$ which depends on the position $q$ of the
tool edge. In [12] the function $f$ is given by
$$f\left(\widt q,\dot q,\beta\right)=-\f{2\pi k_1}{8\beta\Om m}\left[\left(\dot q_1-
\widt q_1\right)+\f{5}{\beta}\left(q_1-\widt
q_1\right)^3\right],\eqno(3.20)$$ where $k_1=\ds\f{3}{4}a
w\beta^{-\f{1}{4}}$ is the parameter depending on further
technological parameters and considered to be constant, $w$ is the
width of the chip, $\Om$ is the constant angular velocity of
rotating work--piece (or tool), the delay $\tau=\ds\f{2\pi}{\Om}$
is the time period of one revolution and $\widt
q{}^1(t)=q^1(t-\tau)$. For $L\left(q,\dot q\right)=\ds\f{1}{2}\dot
q{}^2-\ds\f{1}{2}\al q^2$ and the force filed with delay given by
$$Y\left(q,\dot q\right)=-2k\al \dot q+f\left(q,\dot q,\beta\right),\eqno(3.21)$$
the equation (3.19) is the Euler--Lagrange equation with force with delay. The model was analysed in [12]. The
force field given by (3.21) is not dissipative.

\section*{\normalsize\bf 4. Equivariant dissipation with delay}

\hspace{0.6cm} In this section we study Lagrangian systems that
are invariant under a group action and we shall add to them, in
the sense of the preceding sections, dissipative fields with delay
that are equivariant. This invariance property will yield
dissipative mechanisms with delay that preserve the basic
conserved quantities, yet dissipate enery, as we shall see.

Let $G$ be a Lie group acting on the configuration manifold $Q$ and assume that the lifted action leaves the
Lagrangian $L$ invariant. In this case, the fibre derivative $\ff L:TQ\to T^*Q$ is equivariant with respect to
this action on $TQ$ and the dual action on $T^*Q$. The action $A$, the energy $E$ and the Lagrangian $2$--form
$\Om_L$ are  all invariant under the action of $G$ on $TQ$. Let $Z$ be the Lagrangian vector field for the
Lagrangian $L$, which we assume to be regular. Because of regularity, the vector field $Z$ is also invariant under
$G$. If the action is free and proper, so that $(TQ)/G$ is a manifold, then the vector field and its flow $F_t$
drop to a vector field $Z^G$ and a flow $F_t^G$ on $(TQ)/G$, [2].

Let $J:TQ\to {\bf g}^*$ be the momentum map associated with the
$G$--action, given by:
$$J(v_q)\cdot\xi=<\ff L(v_q),\xi_Q(q)>\eqno(4.1)$$
for $v_q\in T_q Q$ and for $\xi\in{\bf g}$, where $\xi_Q$ denotes
the infinitesimal generator for the action on $Q$. The
infinitesimal generator for the action on the tangent bundle will
be likewise denoted  by $\xi_{TQ}$ and for later use, we note the
relation $T\pi\circ\xi_{TQ}=\xi_{Q}\circ\pi$. If $v(t)$ denotes an
integral curve of the vector field with an equivariant dissipation
term with delay $Y$ added as in the preceding section and if
$J^\xi(v)=<J(v),\xi>$ is the $\xi$--component of the momentum
mapping, then we have:
$$\f{d}{dt}J^\xi(v(t))=dJ^\xi(v(t))\cdot Z(v(t))+dJ^\xi(v(t))Y\left(v(t),\widt v(t)\right).\eqno(4.2)$$
The first term vanishes by conservation of the momentum map for
the Lagrangian vector field $Z$. From (2.8) and the definition of
the momentum map, we get:
$$\begin{array}{lll}
\vspace{0.1cm} dJ^\xi(v)\cdot Y\left(v,\widt v\right)&=&
(i_{\xi_{TQ}}\Om_L)(Y)\left(v,\widt v\right)=-(i_Y\Om_L)
(\xi_{TQ})\left(v,\widt v\right)\\
\vspace{0.1cm} &=& <F^Y\left(v,\widt v\right), T_v T\pi
(\xi_{TQ}(v))>=\\
&=&<F^Y\left(v,\widt
v\right),\xi_Q(\pi(v))>\end{array}\eqno(4.3)$$ and therefore
$$
\f{d}{dt} J^\xi(v(t))=<F^Y\left(v(t),\widt v(t)\right),\xi_Q T\pi
(v(t))>.\eqno(4.4)$$ We summarize this discussion as follows.

{\bf Proposition 4.1.} {\it The momentum map $J:TQ\to{\bf g}^*$ is
conserved under the flow of a $G$--invariant dissipative vector
field with delay $Z+Y$ if and only if $<F^Y,\xi_Q\circ \tau>=0$
for all Lie algebra elements $\xi\in {\bf g}$.}

In this paper we shall consider dissipative vector fields with
delay for which the flow drops to the reduced spaces. Thus a first
requirement is that $Y$ be a vertical $G$--invariant vector field
on  $TQ\times TQ$. A second requirement is that all integral
curves $v(t)$ of $Z+Y$ preserve the sets $J^{-1}{\calO}$, where
${\calO}$ is an arbitrary coadjoint orbit in ${\bf g}^*$. Under
these hypotheses the vector field $Z+Y$ induces a vector field
$Z^G+Y^G$ on $(TQ\times TQ)/G\times G$ that preserves the
symplectic leaves of this Poisson manifold, namely all reduced
spaces $J^{-1}{\calO}/G$.

The condition that $v(t)\in J^{-1}{\calO}$ is equivalent to
$J(v(t))\in {\calO}$, i.e., to the existence of an element
$\eta(t)=\eta(v(t))\in {\bf g}$ such that $dJ(v(t))/dt=
ad^*_{\eta(t)}J(v(t))$ or
$$\f{d J^\xi (v(t))}{dt}=J^{[\eta(t),\xi]}(v(t))\eqno(4.5)$$
for all $\xi \in{\bf g}$. In view of (4.4), we get the following

{\bf Corollary 4.2.} {\it The integral curves of the vector field
$Z+Y$, for $Y$ a vertical $G$--invariant vector field on $TQ\times
TQ$ and $Z$ the Lagrangian vector field of a $G$--invariant
Lagrangian function $L:TQ\to \rr$, preserve  the inverse images of
the coadjoint orbits in ${\bf g}^*$ by the momentum map $J$ if and
only if for each $\left(v,\widt v\right)\in TQ\times TQ$ there is
some $\eta\left(v,\widt v\right)\in{\bf g}$ such that
$$<F^Y\left(v,\widt v\right), (\xi_Q\circ T\pi)(v)>=J^{\left[\eta\left(v,\widt v\right),
\xi\right]}(v)\eqno(4.6)$$ for all  $\xi\in {\bf g}$. As before,
$F^Y$ denotes the force field with delay induced by $Y$.}

We shall see in Section 5 how to construct such a force field with
delay in the case $Q=G$.

\section*{\normalsize\bf 5. Dissipation with delay for Euler--Poincar\'e and Lie--Poisson\break equations}

\hspace{0.6cm} A key step in the reduction of the Euler--Lagrange
equations from the tangent bundle $TG$ of a Lie group $G$ to its
Lie algebra ${\bf g}$ is to understand how to drop the variational
principle to the quotient space. The formulation of the
Euler--Poincar\'e equations and the reduced variational principle
is given by

{\bf Theorem 5.1.} {\it Let $G$ be a Lie group and $L:TG\to \rr$ a
left invariant Lagrangian. Let $l:{\bf g}\to\rr$ be its
restriction to the tangent space at the identity. For a curve
$g(t)\in G$, let $\xi(t)=T_{g(t)}L_{g(t)^{-1}}\dot g(t)$. Then the
followings are equivalent:

{\rm {i)}}  $g(t)$ satisfies the Euler--Lagrange equations for $L$
on $G$.

{\rm {ii)}} The variational principle
$$\de\int_a^b L\left(g(t),\dot g(t)\right)dt=0\eqno(5.1)$$
holds for variations with fixed endpoints.

{\rm iii)} The Euler--Poincar\'e equations hold,
$$\f{d}{dt}\left(\f{\p l}{\p\xi}\right)=ad^*_\xi\left(\f{\p l}{\p\xi}\right).\eqno(5.2)$$

{\rm iv)} The variational principle
$$\de\int_a^b l(\xi(t))dt=0\eqno(5.3)$$
holds on ${\bf g}$, using variations of the form
$\de\xi=\dot\eta+[\xi,\eta]$ where $\eta$ vanishes at the
endpoints.}

In coordinates, the Euler--Poincar\'e equations are
$$\f{d}{dt}\left(\f{\p l}{\p\xi^d}\right)=C_{ad}^b\f{\p L}{\p \xi^b}\xi^a\eqno(5.4)$$
where $C_{ad}^b$ are the structure constants of ${\bf g}$ relative
to a given basis $\xi^a$ are the components of $\xi$ relative to
this basis.

Since the Euler--Lagrange and the Hamilton equations on $TQ$ and
$T^*Q$ are equivalent if the fibre derivative of $L$ is a
diffeomorphism from $TQ$ to $T^*Q$, it follows that the
Lie--Poisson and the Euler--Poincar\'e equations are also
equivalent under similar hypotheses. To see this directly, we make
the following Legendre transformation from ${\bf g}$ to ${\bf
g}^*$:
$$\mu=\f{\p l}{\p\xi},\quad h(\mu)=<\mu,\xi>-l(\xi)\eqno(5.5)$$
and assume that $\xi\to\mu$ is a diffeomorphism. Note that
$\ds\f{\p h}{\p\mu}=\xi$ and so it is now clear that the
Euler--Poincar\'e equations are equivalent to the Lie--Poisson
equations on ${\bf g}^*$, namely
$$\f{d\mu}{dt}=ad^*_{\f{\p h}{\p\mu}}\mu.\eqno(5.6)$$
Now we are ready to synthesize our discussions on forces with
delay and on the Euler--Poincar\'e equations and to transfer this
forcing with delay to the Lie--Poisson equations by means of the
Legendre transform. We begin with a formulation of the
Lagrange--d'{}Alembert principle with delay.

{\bf Theorem 5.2.} {\it Let $G$ be a Lie group, $L:TG\to\rr$ a
left invariant Lagrangian, and $F:TG\times TG\to T^*G$ a force
field with delay equivariant relative to the canonical left
actions of $G$ on $TG\times TG$ and $T^*G$ respectively. Let
$l:{\bf g}\to\rr$ and $f:{\bf g}\times{\bf g}\to{\bf g}^*$ be the
restriction of $L$ and $F$ to $T_eG={\bf g}$. For a curve $g(t)\in
G$, let $\xi(t)=T_{g(t)}L_{g(t)^{-1}}\dot g(t)$. Then the
followings are equivalent:

{\rm i)} $g(t)$ satisfies the Euler--Lagrange equations with
forcing with delay for $L$ on $G$.

{\rm ii)} The integral Lagrange--d'{}Alembert principle with delay
$$\de\int_a^b L(g(t),\dot g(t))dt=\int_a^b F\left((g(t),\dot g(t)),\left(\widt g(t),
\dot{\widt g}(t)\right)\right)\cdot \de g(t)dt\eqno(5.7)$$ holds
for all variations $\de g(t)$ with fixed endpoints.

{\rm iii)} The Euler--Poincar\'e equations with forcing with delay
are valid
$$\f{d}{dt}\left(\f{\p l}{\p\xi}\right)-ad_\xi^*\f{\p l}{\p\xi}=f\left(\xi,\widt\xi\right).\eqno(5.8)$$

{\rm iv)} The variational principle
$$\de\int_a^b l(\xi(t))dt=\int_a^b f\left(\xi(t),\widt \xi(t)\right)\cdot \de\xi(t)dt\eqno(5.9)$$
holds on ${\bf g}$ using variations of the form
$\de\xi=\dot\eta+[\xi,\eta]$ where $\eta$ vanishes at the
endpoints.}

{\bf Proof.} We have already seen that i) and ii) are equivalent for any confi\-gu\-ra\-tion manifold $Q$ in
Section 2. Next we prove that ii) and iv) are equivalent. First, note that $l:{\bf g}\to \rr$ and $f:{\bf
g}\times{\bf g}\to{\bf g}^*$ determine uniquely a function $L:TG\to\rr$, and a function $F:TG\times TG\to TG$ by
left translation of the argument and conversely. Thus, the equivalence of ii) and iv) comes down to proving that
all variations $\de g(t)\in TG$ of $g(t)$ with fixed endpoints induce and are induced by variations $\de\xi(t)$ of
$\xi(t)$ of the form $\de\xi=\dot\eta+[\xi,\eta]$ where $\eta(t)$ vanishes at the endpoints. But this is precisely
the matter of the proposition 5.1 in [2].

The Euler--Poincar\'e equations with forcing with delay have the
following expression in local coordinates:
$$\f{d}{dt}\left(\f{\p l}{\p\xi^a}\right)-C_{ba}^d\f{\p l}{\p \xi^d}=f_a\left(\xi,\widt\xi\right)\eqno(5.10)$$
where $C_{ba}^d$ are the structure constants of the Lie algebra
${\bf g}$.

The condition that the integral curves of the dissipative vector
field with delay preserve the inverse images of coadjoint orbits
by momentum map and hence the integral curves of (5.8) preserve
the coadjoint orbits of ${\bf g}^*$ is given by (4.4), Section 4.
Since $\xi_G(g)=T_eR_g(\xi)$ and $J(v_g)=T_e^*R_g\ff L(v_g)$, we
get $$<F,\xi_G\circ\pi>\left(v_g,\widt v_{\tilde
g}\right)=<F\left(v_g,\widt v_{\tilde g}\right),
T_eR_g(\xi)>=T_e^*R_g F\left(v_g,\widt v_{\tilde
g}\right)\cdot\xi$$ and
$$\begin{array}{lll}
\vspace{0.1cm}
 J^{\left[\eta\left(v_g,\widt v_{\tilde
g}\right),\xi\right]}(v_g)&=&T_e^*R_g\ff L(v_g)\cdot
\left[\eta\left[v_g,\widt v_{\tilde g}\right),\xi\right]=\\
&=& (ad^*_{\eta(v_g,\widt v_{\tilde g})}\circ T_e^* R_g\circ\ff
L)(v_g)\cdot\xi.\end{array}$$ Since $F$ and $\ff L$ are
equivariant,
$$T_e^* R_g F(v_g,\widt v_{\tilde g})=Ad^*_{g^{-1}} F\left(T_gL_{g^{-1}}v_g, T_{\tilde g}
L_{\tilde g{}^{-1}}\widt v_{\tilde g}\right)$$ and
$$\left(ad^*_{\eta\left(v_g,\widt v_{\tilde g}\right)}\circ T_e^* R_g\circ\ff L\right)(v_g)=
\left(ad^*_{\eta\left(v_g,\widt v_{\tilde g}\right)}\circ
Ad^*_{g^{-1}}\circ\ff l\right) (T_g L_{g^{-1}}v_g).$$ Because
$Ad_{g^{-1}}\circ ad_{\eta\left(v_g,\widt v_{\tilde
g}\right)}=ad_{Adg^{-1}\eta\left(v_g,\widt v_{\tilde
g}\right)}\circ Ad_{g^{-1}}$ we get
$$J^{\left[\eta\left(v_g,\widt v_g\right),\xi\right]}(v_g)=Ad^*_{g^{-1}}
\circ ad^*_{Ad_{g^{-1}}\eta\left(v_g,\widt v_{\tilde
g}\right)}\circ\ff L)(T_gL_{g^{-1}}v_g)$$ and the identity (4.4),
Section 4 thus becomes
$$F\left(T_g L_{g^{-1}}v_g, T_{\widt g}L_{\widt g{}^{-1}}\widt v_{\tilde g}\right)=
\left(ad^*_{Ad_{g^{-1}}\eta\left(v_g,\widt v_{\tilde
g}\right)}\circ\ff L\right)\left(T_g L_{g^{-1}}v_g\right).$$
Letting $\xi=T_g L_{g^{-1}}v_g$, $\widt\xi =T_{\widt g}L_{\widt g
{}^{-1}}\widt v_{\tilde g}$, this becomes
$$f\left(\xi,\widt \xi\right)=ad^*_{Ad_{g^{-1}}\eta\left(v_g,\widt v_{\tilde g}\right)}\f{\p l}{\p\xi}(\xi).\eqno(5.11)$$

The left hand side is independent of $g$ and thus the right hand
side must be also $g$--independent. Thus taking $g=\widt g=e$ the
criterion (4.4), Section 4 becomes: for every $\xi,\widt\xi\in{\bf
g}$ there is some $\eta\left(\xi,\widt\xi\right)\in{\bf g}$ such
that
$$f\left(\xi,\widt\xi\right)=ad^*_{\eta\left(\xi,\widt\xi\right)}\f{\p l}{\p\xi}(\xi).\eqno(5.12)$$
In other words, the force field with delay $f$ (and hence $F$) is
completely determined by an arbitrary map $\eta:{\bf g}\times{\bf
g}\to{\bf g}$ via formula (5.11) and we conclude the following

{\bf Corollary 5.3.} {\it The solutions of the Euler--Poincar\'e
equations with for\-cing with delay $(5.8)$ perserve the coadjoint
orbits of ${\bf g}^*$ provided the force field with delay $f$ is
given by $(5.12)$ for some smooth map $\eta:{\bf g}\times{\bf
g}\to{\bf g}$.}

Transforming the Euler--Poincar\'e equations with forcing with
delay by means of the Legendre transformation (5.5), the equations
(5.8) with the force field with delay (5.12) become $$\f{
d\mu}{dt}-ad^*_{\frac{\p
l}{\p\mu}}\mu=-ad^*_{\eta\left(\mu,\widt\mu\right)}\mu\eqno(5.13)$$
where $\eta:{\bf g}^*\times{\bf g}^*\to{\bf g}$. The requirement
on the map $\eta$ is that the right hand side of (5.13) be a
gradient relative to a certain metric on the orbit.

To generalize the metric defined by Killing form [3] to coadjoint orbits of the dual ${\bf g}^*$ of a general Lie
algebra ${\bf g}$, we introduce a symmetric positive definite bilinear form.

Let a symmetric positive definite bilinear form $\widt\Gamma:{\bf
g}^*\times{\bf g}^*\to\rr$ and denote by $\Gamma:{\bf g}^*\to{\bf
g}$ the induced map given by
$\widt\Gamma(\al,\beta)=<\beta,\Gamma\al>$ for all
$\al,\beta\in{\bf g}^*$, where $<,>:{\bf g}^*\times{\bf g}\to\rr$
is the pairing between ${\bf g}^*$  and ${\bf g}$. Symmetry of
$\widt\Gamma$ is equivalent to symmetry of $\Gamma$, i.e.
$\Gamma^*=\Gamma$. We introduce the following new inner product on
${\bf g}$:
$$<\xi,\eta>_{\Gamma^{-1}}\,=\,<\Gamma^{-1}\eta,\xi>\eqno(5.14)$$
for all $\xi,\eta\in{\bf g}$ and call it the $\Gamma^{-1}$--{\it
inner product}.

Let ${\bf g}_\mu=\{\xi\in{\bf g}~|~ad^*_{\xi}\mu=0\}$ denote the
coadjoint isotropy subalgebra of $\mu\in{\bf g}$ and denote by
${\bf g}^\mu$ the orthogonal complement of ${\bf g}_\mu$ relative
to the $\Gamma^{-1}$--inner product. For an element $\xi\in{\bf
g}$ we denote by $\xi_\mu$ and $\xi^\mu$ the components of $\xi$
in the orthogonal direct sum decomposition ${\bf g}={\bf
g}_\mu\oplus {\bf g}^\mu$. For all $g\in G$, $\xi\in{\bf g}$ and
$Ad^*_g:{\bf g}^*\to{\bf g}^*$, $ad^*_\xi:{\bf g}^*\to{\bf g}^*$,
$Ad_g:{\bf g}\to{\bf g}$ we have
$$\begin{array}{l}
\vspace{0.1cm}
\Gamma(Ad_g^*\mu)=Ad_g\Gamma\mu,\\
ad_\xi^*(Ad_g^*\mu)=Ad_g^*(ad^*_{Ad_g\xi}(Ad^*_{g^{-1}}\mu)),\quad
\mu\in {\bf g}^*.\end{array}\eqno(5.15)$$

Let $\calO_{\mu_0}$ be the coadjoint orbit through $\mu_0\in {\bf
g}^*$ and $\mu,\widt\mu\in\calO_{\mu_0}$. There exists $g\in G$
such that $\mu=Ad^*_{g^{-1}}\widt\mu$ and ${\bf
g}_\mu=Ad_{g^{-1}}{\bf g}_{\widt\mu}$, ${\bf
g}^\mu=Ad_{g^{-1}}{\bf g}^{\widt\mu}$.

Let $C$ be a positive Casimir function on ${\bf g}^*$ and let
$\om$ be the coadjoint orbit symplectic structure defined by
$$\om(\mu)(\xi_*(\mu),\eta_*(\mu))=-<\mu,[\xi,\eta]>\eqno(5.16)$$
for all $\mu\in\calO_{\mu_0}$, $\xi,\eta\in{\bf g}$ and
$\xi_*(\mu),\eta_*(\mu)\in T\calO_{\mu_0}$. If
$\mu,\widt\mu\in\calO_{\mu_0}$ then $ad^*_\xi\mu\in
T_\mu\calO_{\mu_0}$, $ad^*_\xi\widt\mu\in
T_{\widt\mu}\calO_{\mu_0}$. Since $\mu=Ad^*_{g^{-1}}\widt\mu$ it
results
$$ad^*_\xi\mu=Ad^*_{g^{-1}}\left(ad^*_{Ad_g\xi}\widt\mu\right).\eqno(5.17)$$

We define the $(C,\Gamma^{-1})$--{\it normal metric} on
$\calO_{\mu_0}$  {\it with respect to $\mu$ and $\widt\mu$} by
$$<ad^*_\xi\widt\mu, ad_\eta^*\widt\mu>_N(\mu)=C(\widt\mu)<\Gamma^{-1}\eta^{\widt\mu},
\xi^{\widt\mu}>-$$
$$-\f{1}{C(\mu)}\om(\Gamma\xi)\left((\Gamma\mu)_*,
(\Gamma\widt\mu)_*\right)\cdot\om(\Gamma\eta)\left((\Gamma\mu)_*,(\Gamma\widt\mu)_*\right)\eqno(5.18)$$
for all $\xi,\eta\in {\bf g}$.

We shall regard $C$ and $\Gamma$ as fixed in the following discussion and just refer to this metric as the normal
metric. If $\mu=\widt\mu$ then the normal metric is given in [2].

Let $k:{\bf g}^*\to\rr$ be a smooth function. We shall compute the
gradient vector of $k_{|\calO_{\mu_0}}$ with respect to the normal
metric. For this purpose we denote by $\ds\f{\p k}{\p \mu}\in {\bf
g}$ the derivative of $k$ at $\mu$ and by $\grad k(\mu)$ the
gradient of $k_{|\calO_{\mu_0}}$. Since $\grad k(\mu)\in
T_\mu\calO_{\mu_0}$ we can write $\grad k(\mu)=ad^*_\eta\mu$ for
some $\eta\in{\bf g}$. Since $\xi_\mu$ and $\eta^\mu$ are
orthogonal in the $\Gamma^{-1}$--inner product, we get
$$-<ad^*_{\frac{\p k}{\p\mu}}\mu,\xi>(\mu)=\left<\mu,\left[\xi,\f{\p k}{\p\mu}\right]\right>(\mu)=\left< ad^*_\xi\mu,\f{\p k}{\p\mu}\right>(\mu)=$$
$$=<\grad k(\mu), ad^*_\xi\mu>_N(\mu)=<ad^*_\eta\mu, ad^*_\xi\mu>_N(\mu)=C(\widt\mu)<\Gamma^{-1}\xi^{\widt\mu},\eta^{\widt\mu}>-$$
$$-\f{1}{C(\mu)}<\Gamma\xi,[\Gamma\mu,\Gamma\widt\mu]>\cdot <\Gamma\eta,[\Gamma\mu,\Gamma\widt\mu]>=C(\widt \mu)<\Gamma^{-1}(\xi^{\widt\mu}+\xi_{\widt\mu}),\eta^{\widt\mu}>-$$
$$-\f{1}{C(\mu)}<\Gamma^{-1}[\Gamma\mu,\Gamma\widt\mu],\xi>\cdot<\Gamma^{-1}[\Gamma\mu,\Gamma\widt\mu],\eta>=C(\widt\mu)<\Gamma^{-1}\eta^{\widt\mu},\xi>-$$
$$-\f{1}{C(\mu)}<\Gamma^{-1}[\Gamma\mu,\Gamma\widt\mu],\eta>\cdot<\Gamma^{-1}[\Gamma\mu,\Gamma\widt\mu],\xi>$$
for any $\xi\in{\bf g}$. Therefore
$$C(\widt\mu)\Gamma^{-1}\eta^{\widt\mu}-\f{1}{C(\mu)}\Gamma^{-1}[\Gamma\mu,
\Gamma\widt\mu]^\mu<\Gamma^{-1}[\Gamma\mu,\Gamma\widt\mu],\eta>=-ad^*_\eta\mu.\eqno(5.19)$$

Because $\eta^{\widt\mu}=Ad_g\eta^\mu$ it results
$$C(Ad_g^*\mu)Ad_g\eta^\mu-\f{1}{C(\mu)}\left<\Gamma\left(\f{\p k}{\p\mu}\right),[\Gamma\mu,\Gamma\widt\mu]\right>=-
\Gamma\left(ad^*_{\frac{\p k}{\p\mu}}\mu\right)\eqno(5.20)$$
and
$$\eta^\mu=-\f{1}{C(\widt\mu)}Ad_{g^{-1}}\Gamma\left(ad^*_{\f{\p
k}{\p\mu}}\mu\right)+\f{1}{C(\mu)C(\widt\mu)}\left<\Gamma\left(\f{\p
k}{\p\mu}\right),[\Gamma\mu,\Gamma\widt\mu]\right>\cdot[\Gamma\mu,\Gamma\widt\mu]^\mu$$
or
$$\eta^\mu=-\f{1}{C(\widt\mu)}\Gamma\left(ad^*_{\f{\p k}{\p\widt\mu}}\widt\mu\right)+
\f{1}{C(\mu)C(\widt\mu)}\left<\Gamma\left(\f{\p k}{\p\mu}\right),[\Gamma\mu,\Gamma\widt\mu]\right>
\cdot [\Gamma\mu,\Gamma\widt\mu]^\mu.\eqno(5.21)$$ Thus
$$\grad k(\mu)=ad^*_{\eta^\mu}\mu=-\f{1}{C(\widt\mu)}\Gamma\left(ad^*_{\f{\p k}{\p\widt\mu}}\widt\mu
\right)+$$ $$+\f{1}{C(\mu)C(\widt\mu)}\left<\Gamma\left(\f{\p
k}{\p\mu}\right),[\Gamma\mu,\Gamma\widt\mu]\right>\cdot[\Gamma\mu,\Gamma\widt\mu]^\mu\eqno(5.22)$$
and the equation of the gradient vector field in
$\mu\in\calO_{\mu_0}$ relative to the normal metric on
$\calO_{\mu_0}$ is
$$\f{d\mu}{dt}=-\f{1}{C(\widt\mu)}ad^*_{\Gamma\left(\f{\p k}{\p\widt\mu}\right)}\widt\mu+\f{1}{C(\mu)C(\widt\mu)}\left<\Gamma\left(\f{\p k}{\p\mu}\right),[\Gamma\mu,\Gamma\widt\mu]\right>\cdot [\Gamma\mu,\Gamma\widt\mu]^\mu.\eqno(5.23)$$

Therefore in (5.23) we put $\eta(\mu,\widt\mu)=-\eta^{\mu}$ and
the Lie--Poisson equations with delay forcing (5.6) become
$$\f{d\mu}{dt}\!=\!\f{1}{C(\widt\mu)}ad^*_{\Gamma\Bigl(ad^*_{\f{\p k}{\p\widt\mu}}
\widt\mu\Bigr)}\mu
\!-\!\f{1}{C(\mu)C(\widt\mu)}<\Gamma\!\!\left(\!\f{\p
k}{\p\mu}\!\right),[\Gamma\mu,\Gamma\widt\mu]>\!\cdot
ad^*_{[\Gamma\mu,\Gamma\widt\mu]^\mu}\mu.\eqno(5.24)$$

If ${\bf g}$ is a compact algebra with the bi-invariant inner
product $<\cdot,\cdot>$ on ${\bf g}$ and ${\bf g}$ is also
semisimple, then we could let $<\cdot,\cdot>=-k(\cdot,\cdot)$
where $-k(\cdot,\cdot)$ is the Killing form. In these conditions
the inner product identifies ${\bf g}$ with its dual ${\bf g}^*$,
coadjoint orbits with adjoint orbits so that
$ad^*_{\xi}\mu=[\mu,\xi]$ and $\ds\f{\p k}{\p\mu}=\nabla k(\mu)$,
where $\nabla k(\mu)$ is the gradient of $k$ on ${\bf g}$ at $\mu$
relative to the bi-invariant inner product $<\cdot,\cdot>$. The
formula for the gradient vector field on the adjoint orbit
$\calO_{\mu_0}$ relative to $\mu$ and $\widt\mu$ becomes
$$\f{d\mu}{dt}=-\f{1}{C(\widt\mu)}\left[\mu,\Gamma[\widt\mu,\nabla k(\widt\mu)]\right]-
\f{1}{C(\mu)C(\widt\mu)}\left<\Gamma\left(\f{\p k}{\p\mu}\right),
[\Gamma\mu,\Gamma\widt\mu]\right>\cdot[\mu,[\Gamma\mu,\Gamma\widt\mu]^\mu],\eqno(5.25)$$
where $\Gamma:{\bf g}\to{\bf g}$ defines the symmetric positive
definite bilinear form $(\xi,\eta)\mapsto <\Gamma\xi,\eta>$. Thus
in this case the Lie--Poisson equations with delay forcing become
$$\f{d\mu}{dt}=-[\nabla
h(\mu),\mu]+\f{1}{C(\widt\mu)}[\mu,\Gamma\left[\widt\mu,\nabla
k(\widt\mu)]\right]-$$
$$-\f{1}{C(\mu)C(\widt\mu)}\left<\Gamma(\nabla k(\mu)),[\Gamma\mu,\Gamma\widt\mu]\right>
\cdot[\mu,[\Gamma\mu, \Gamma\widt\mu]^\mu].\eqno(5.26)$$ Taking
$C(\mu)=1$ and $\Gamma$ to be the identity, the dissipative term
with delay in (5.26) is the Brockett double bracket. The condition
that the delay forcing term be dissipative is $\ds\f{dh}{dt}<0$
and this imposes some conditions on the choice of the function
$k:{\bf g}^*\to\rr$.

\section*{\normalsize\bf 6. Free rigid body with delay}

\hspace{0.6cm} Let $G=SO(3)$ and the usual identification
$(so(3),[\cdot,\cdot]\cong (\rr^3, \times)$ implies $so(3)^*\cong
\rr^3$ via the natural pairing given by the Euclidean inner
product. Consider $\calO(M_0)$ the coadjoint orbit through
$M_0\in\rr^3$. The infinitesimal generator of the coadjoint action
is given by $\xi_{so(3)^*}(M)=\xi\times M$ for $M\in{\calO}(M_0)$
and $\xi\in so(3)$. For $M,\widt M\in{\calO}(M_0)$ let the tangent
vectors $\xi_{\rr^3}(M)=\xi\times M\in T_M{\calO}(M_0)$ and
$\xi_{\rr^3}(\widt M)=\xi\times\widt M\in T_{\widt
M}{\calO}(M_0)$; the coadjoint orbit symplectic structure becomes
$$\om\left(\xi_{\rr^3}(M),\xi_{\rr^3}(\widt M)\right)=M\cdot \left(\widt M\times\xi\right).\eqno(6.1)$$
The normal metric on ${\calO}(M_0)$ with respect to $M$ and $\widt
M$, given by (5.18), is
$$\left<\widt M\times \xi,\widt M\times\eta\right>_N(M)=\f{1}{c^4}\left<\widt M\times\left(
\widt M\times\xi\right), \widt M\times\left(\widt
M\times\eta\right)\right>-$$ $$-\f{1}{c^4}\left(M\cdot\left(\widt
M\times\xi\right)\right)\cdot \left(M\cdot\left(\widt
M\times\eta\right)\right)=$$
$$=\left<\widt M\times\xi,\widt M\times\eta\right>_N\left(\widt M\right)-\f{1}{c^4}
\left(M\cdot\left(\widt
M\times\xi\right)\right)\cdot\left(M\cdot\left(\widt M\times
\eta\right)\right)\eqno(6.2)$$ with $c=C(\widt M)=C(M)$.

The normal metric at $\widt M$ on two tangent vectors $\widt
M\times\xi$, $\widt M\times\eta$ to the sphere of radius $c$ is
given by
$$\left<\widt M\times\xi,\widt M\times\eta\right>_N\left(\widt M\right)=\f{1}{c^4}\left<\widt M\times\left(\widt M\times\xi\right),\widt M\times\left(\widt M\times\eta\right)\right>\eqno(6.3)$$
where the inner product of the right hand side is the standard
inner product in $\rr^3$.

The normal metric at $M$ with $M\not=\widt M$, on two tangent
vectors $\widt M\times\xi$, $\widt M\times\eta$ to the sphere of
radius $c$ at $\widt M$ is given by
$$\left<\widt M\times\xi,\widt M\times\eta\right>_N(M)=\f{1}{c^4}\left<M\times\left(\widt M\times\xi\right), M\times\left(\widt M\times\eta\right)\right>\eqno(6.4)$$
where the inner product of the right hand side is the standard
inner product in $\rr^3$. From (6.3) and (6.4) it results (6.2).

From (5.26) where $k(M)=h(M)=\ds\f{1}{2}\|M\|^2$ it results the Lie--Poisson equation for the rigid body with
delay:
$$\dot M(t)=M\times\Om+\f{\al}{c^2}M\times\left(\widt M\times\widt\Om\right).\eqno(6.5)$$

We shall discuss the stability of the equilibrium states for a
free rigid body with delay.

Let the free rigid body with delay given by the equation$$\dot
M=M\times\Om+\al M\times\left(\widt
M\times\widt\Om\right),\eqno(6.6)$$ where $M=I\Om=\left(I_1 x(t),
I_2y(t), I_3 z(t)\right)^T$, $\Om=\left(x(t),y(t), z(t)\right)^T$,
$\widt M=I\widt\Om=\left(I_1x(t-\tau), I_2 y(t-\tau), I_3
z(t-\tau)\right)^T$, $\widt\Om=\left(x(t-\tau), y(t-\tau),
z(t-\tau)\right)^T$, $I_1>0$, $I_2>0$, $I_3>0$, $\al$ a constant
and $\tau\ge 0$.

It is not hard to see that the equilibrium states of our system
are $\Om_1=(m/I_1,0,0)^T$, $\Om_2=(0, m/I_2,0)^T$, $\Om_3=(0,0,
m/I_3)$, $m\in\rr^*$.

{\bf Proposition 6.1.} {\it The equilibrium state $\Om_1$ has the
following behavior:

{\rm (i)} The corresponding linear system is given by
$$\de\dot\Om=A\de\Om+\al G\de\widt\Om\eqno(6.7)$$
where
$$A=\left(\begin{array}{ccc}
\vspace{0.2cm}
0 & 0 & 0\\
\vspace{0.2cm}
0 & 0 & \f{I_3-I_1}{I_1I_2}m\\
0 & \f{I_1-I_2}{I_1I_3}m & 0\end{array}\right),\quad
G=\left(\begin{array}{ccc}
\vspace{0.2cm}
0 & 0 & 0\\
\vspace{0.2cm}
0 & \f{I_2-I_1}{I_1I_2}m^2 & 0\\
0 & 0 & \f{I_3-I_1}{I_1I_3}m^2\end{array}\right);\eqno(6.8)$$

{\rm (ii)} The characteristic equation is
$$\lam\left[\lam^2-\f{\al m^2}{I_1}\left(\f{I_2-I_1}{I_2}+\f{I_3-I_1}{I_3}\right)\lam e^{-\tau\lam}+\f{\al^2m^4}{I_1^2I_2I_3}(I_2-I_1)(I_3-I_1)e^{-2\tau\lam}-\right.$$
$$-\left.\f{(I_1-I_2)(I_3-I_1)}{I_1^2I_2I_3}m^2\right]=0;\eqno(6.9)$$

{\rm (iii)} On the tangent space at $\Om_1$ to the sphere of
radius $m^2$ the linear operator given by the linearized vector
field has the characteristic equation
$$\lam^2-\f{\al m^2}{I_1}\left(\f{I_2-I_1}{I_2}+\f{I_3-I_1}{I_3}\right)\lam e^{-\tau\lam}+\f{\al^2m^4}{I_1^2I_2I_3}(I_2-I_1)(I_3-I_1)e^{-2\tau\lam}-$$
$$-\f{(I_1-I_2)(I_3-I_1)}{I_1^2I_2I_3}m^2=0;\eqno(6.10)$$

{\rm (iv)} If $I_1>I_2$, $I_1>I_3$ for $0\le\tau <\tau_c$, where
$$\tau_c=\f{I_1\left[I_3(I_1-I_2)+I_2(I_1-I_3)\right]}{3|\al|m^2(I_1-I_2)(I_1-I_3)},\eqno(6.11)$$
then the equilibrium state $\Om_1$ is asymptotically stable.}

{\bf Proof.} (i), (ii), (iii) result from the definitions of the liniarized and the characteristic equation by
calculus. For (iv) consider $I_1>I_2$, $I_1>I_3$. If $\tau=0$ then the characteristic equation (6.9) has
eigenvalues with the real parts negative and $\Om_1$ is asymptotically stable. Following [9] and [12] it results
that for $0<\tau\le\tau_c$ the equilibrium state $\Om_1$ remains asymptotically stable.

In the following we study the existence of Hopf bifurcations for
the free rigid body with delay (6.6) by choosing the delay $\tau$
as a bifurcation para\-meter. First we would like to know when the
equation (6.10) has purely imaginary roots $\pm i\om_0$
$(\om_0>0)$ at $\tau=\tau_0$. Note that $\lam=i\om_0$ is a root of
(6.1) if
$$\begin{array}{l}
\vspace{0.1cm}
\om_0^2-c-a\om_0\sin\om_0\tau_0-b\cos 2\om_0\tau_0=0,\\
a\om_0\cos\om_0\tau_0-b\sin 2\om_0\tau_0=0\end{array}\eqno(6.12)$$
with $a,b,c$ given by
$$a=\f{\al m^2}{I_1^2I_2I_3}[I_3(I_1-I_2)+I_2(I_1-I_3)],\quad b=\f{\al^2 m^4}{I_1^2I_2I_3}(I_1-I_2)(I_1-I_3),$$
$$c=\f{m^2}{I_1^2I_2I_3}(I_1-I_2)(I_1-I_3).\eqno(6.13)$$

We deduce the following

{\bf Proposition 6.2.} (i) {\it If $|m|<\ds\f{1}{|\al|}$ then
$\lam=i\om_0$ is a simple root of $(6.10)$ and
$$\om_0=\f{a+\sqrt{a^2-4(b-c)}}{2},\quad
\tau_0=\f{\pi}{2\om_0},\quad \tau_0>\tau_c.\eqno(6.14)$$

{\rm (ii)} If $|m|>\ds\f{1}{|\al|}$ then $\lam=i\om_0$ is a simple
root of $(6.10)$ and}
$$\om_0=\f{-a+\sqrt{a^2-4(b-c)}}{2},\quad
\tau_0=\f{3\pi}{2\om_0},\quad \tau_0>\tau_c.\eqno(6.15)$$

We proceed to calculate $\Rem\left(\ds\f{d\lam}{d\tau}\right)$ at
$\tau=\tau_0$. By differentiating the equation (6.10) implicitly
with respect to $\tau$, we obtain:
$$\f{d\lam}{d\tau}=\f{a\lam^2 e^{-\tau\lam}+2b\lam e^{-2\tau\lam}}{2\lam+a(1-\tau\lam)e^{-\tau\lam}
-2b\tau e^{-2\tau\lam}}.\eqno(6.16)$$

It is then evaluated at $\lam=i\om_0$ and $\tau=\tau_0$ given by
(6.14), (6.15), yielding
$$\Rem\left(\f{d\lam}{d\tau}\right)_{\lam=i\om_0\atop{\tau=\tau_0}}=\f{\om_0(\om_0+a)(a-2b)}{\tau_0(a\om_0
-2b)^2+(\om_0+a)^2}.\eqno(6.17)$$

From the standard Hopf bifurcation theory we have the following
result.

{\bf Proposition 6.3.} If $I_1>I_2$, $I_1>I_3$ and $\om_0,\tau_0$
are given by (6.14), (6.15) with (6.13), then
$\Rem\left(\ds\f{d\lam}{d\tau}\right)_{\lam=i\om_0\atop{\tau=\tau_0}}\not=0$
and a Hopf--type bifurcation occurs at $\Om_1$ when $\tau$ passes
through $\tau_0$.

In the following we obtain some conditions which guarantee that the free rigid body with delay undergoes a Hopf
bifurcation at $\tau=\tau_0$. The method we use is based on the normal form theory and the center manifold theorem
introduced in [7].

With the translation $V=\Om-\Om_1$, $U=M-M_1$ the equation (6.6)
becomes
$$\dot U=IAV+\al IG\widt V+F\left(U, V,\widt U, \widt V\right),\eqno(6.18)$$
where $A,G$ are given by (6.8), $I=\diag(I_1,I_2,I_3)$ and
$$F\left(U,V,\widt U,\widt V\right)=U\times V+\al\left[U\times\left(\widt U\times\Om_1\right)+
U\times\left(M_1\times\widt V\right)+\right.$$
$$+\left.M_1\times\left(\widt U\times\widt V\right)\right]+\al
U\times\left(\widt U\times\widt V\right).\eqno(6.19)$$

From (6.18) it results
$$\dot V=AV+\al G\widt V+N\left(V,\widt V\right),\eqno(6.20)$$
where $N\left(V,\widt V\right)=I^{-1}F\left(I^{-1}V, V,
I^{-1}\widt V,\widt V\right).$

For $\phi\in C^1([-\tau_0, 0],\rr^3)$ we define an operator $\ca$
by
$$\ca\phi(\th)=\left\{\begin{array}{ll}
\vspace{0.2cm}
\ds\f{d\phi}{d\th}, & \th\in [-\tau_0,0)\\
A\phi(0)+\al G\phi(-\tau_0), &
\th=0\end{array}\right.\eqno(6.21)$$ and for $\psi\in
C^1([0,\tau_0],\rr^{3*})$ we define the adjoint operator ${\ca}^*$
of ${\ca}$ by
$$\ca^*\psi(s)=\left\{\begin{array}{ll}
\vspace{0.2cm}
-\ds\f{d\psi}{ds}, & s\in [0,\tau_0)\\
\psi(0)A+\al\psi(\tau_0)G, &
s=\tau_0;\end{array}\right.\eqno(6.22)$$ $\ca$ and $\ca^*$ are
adjoint operators with respect to the bilinear form
$$<\psi,\phi>=\ov\psi(0)\phi(0)-\al\int_{-\tau_0}^0\int_{\xi=0}^\th\ov\psi(\xi-\th)
G\phi(\xi)d\xi d\th,\eqno(6.23)$$ $\phi\in
C^1([-\tau_0,0],\rr^3)$, $\psi\in C^1([0,\tau_0],\rr^{3*})$.

Let now $\lam_1=i\om_0$, $\lam_2=\ov\lam_1=-i\om_0$ be eigenvalues
of $\ca$, where $\om_0$ is given by (6.14). They are also
eigenvalues of $\ca^*$. We can easily obtain that
$$\phi(\th)=(0, v_2, v_3)^T e^{\lam_1\th},\quad \th\in [-\tau_0,0],\eqno(6.24)$$
where $v_2=(I_3-I_1)m$, $v_3=\lam_1I_1I_2-(I_2-I_1)m^2\al
e^{-\lam_1\tau_0}$, is an eigenvector of $\ca$ corresponding to
$\lam_1$ and
$$\psi(s)=(0, w_2, w_3)e^{\lam_1 s},\quad s\in [0,\tau_0],\eqno(6.25)$$
where $w_2=I_2(I_1-I_2)m$, $w_3=(\lam_1I_1I_2-(I_2-I_1)m^2\al
e^{\lam_1\tau_0})I_3$, is an eigenvector of $\ca^*$ corresponding
to $\lam_1$.

From (6.23), (6.24), (6.25) it results
$$\begin{array}{lll}
\vspace{0.2cm} a_{11}&=&<\psi,\phi>=v_2w_2+v_3\ov w_3-\ds\f{\al
m^2}{I_1^2I_2I_3}\cdot\f{\al \tau_0}{\lam_2^2}
(e^{\lam_2\tau_0}+\\
\vspace{0.2cm}
&+&\lam_2 e^{\lam_2\tau_0}-1)[I_3(I_2-I_1)v_2w_2+I_2(I_3-I_1)v_3\ov w_3],\\
\vspace{0.2cm} a_{12}&=&<\psi,\ov\th>=v_2w_2+\ov v_3\ov
w_3-\ds\f{\al m^2}{I_1^2I_2I_3}\cdot\ds
\f{\al}{2\lam_2^2}(2-e^{-\lam_2\tau_0}-\\
\vspace{0.2cm}
&-&e^{\lam_2\tau_0})[I_3(I_2-I_1)v_2w_2+I_2(I_3-I_1)\ov v_3\ov w_3],\\
\vspace{0.2cm} a_{21}&=&<\ov\psi,\phi>=\ov{<\psi,\ov\phi>}=\ov
a_{12},~a_{22}=<\ov\psi,\ov\phi>=\\
&=&\ov{<\psi,\phi>}=\ov a_{11}.\end{array}\eqno(6.26)$$

Let $d=a_{11}\ov a_{11}-a_{12}\ov a_{12}$ and
$b_{11}=\ds\f{a_{11}}{d}$, $b_{12}=-\ds\f{a_{12}}{d}$. The vector
field given by $\widt \psi(s)=b_{11}\psi(s)+b_{12}\psi(s)$, $s\in
[0,\tau_0]$ is an eigenvector of $\ca^*$  satisfying the
relations:
$$<\widt\psi,\th>=1,\quad <\widt\psi,\ov\th>=<\ov{\widt\psi},\th>=0,\quad <\ov{\widt\psi},\ov\th>=1.\eqno(6.27)$$
Then we have
$$\widt\psi(s)=\left(0,\widt w_2,\widt w_3\right)e^{\lam_1s},\quad
s\in [0,\tau_0],\eqno(6.28)$$ where
$$\widt w_2=(b_{11}+b_{12})w_2,\quad \widt w_3=b_{11}w_3+b_{12}\ov
w_3.\eqno(6.29)$$


We shall discuss the existence of a local center manifold around
the equilibrium point of the equation (6.20).

Let $\cb=C^1([-\tau_0,0],\rr^3)$ and $\De$ the vectorial space,
span of the eigenvectors $\phi(\th)$, $\ov{\phi(\th)}$
corresponding to $\lam_1=i\om_0$, $\lam_2=-i\om_0$. For a given
neighborhood $\cv$ of $0\in\cb$, a local center manifold
$W_\loc^c(0)=W^c(0,\cv)$ of the equilibrium point $O(0,0,0)$ of
(6.20) is a $C^1$--submanifold that is a graph over $\cv\cap\De$
in $\cb$, tangent to $\De$ at $O$ and locally invariant under the
flow defined by the equation (6.20). In other words
$$W_\loc^c(0)=\left\{\fii^c\in\cb~|~\fii^c=u\phi+\ov u\ov\phi+w(\fii),~\fii\in\cv\cap\De\right\},\eqno(6.30)$$
where $w:\De\to\cb$ is a $C^1$--mapping with $\fii(0)=0$, $D_\fii
w(0)=0$ and $\left<\ov{\widt\psi}, w\right>=0$. Moreover, every
orbit that begins on $W_\loc^c(0)$ remains in this set as long as
it stays in $\cv$.

The basic result on the existence of the local center manifold for the delay differential equations is given in
[6]. From the definition of the local center manifold it results that
$$W_\loc^c(0)\cap V_1=\left\{\fii^c\in\cb~|~\fii^c=u\phi+\ov u\ov\phi+w(u,\ov u),~u=u_1+iu_2,\right.$$
$$\left.(u_1,u_2)\in V_1\subset\rr^2\right\},\eqno(6.31)$$ where
$w:[-\tau,0]\times\cc^2\to\rr^3$ is given by $w(\th,u,\ov u)=w\left(u\phi(\th)+\ov u\ov\phi(\th)\right)$.
Following [6], $\fii^c\in W_\loc^c(0)\cap V_1$ and the function $w$ is the solution of the partial derivate system
$$\f{\p w}{\p t}\left(\th, u(t),\ov u(t)\right)+g\left(u(t),\ov u(t)\right)\phi(\th)+\ov{g(u(t),\ov u(t))}\cdot \phi(\th)=\f{\p w}{\p\th}\left(\th, u(t),\ov u(t)\right)\eqno(6.32)$$
with
$$\f{\p w}{\p t}\left(0, u(t),\ov u(t)\right)+g\left(u(t),\ov u(t)\right)\phi(0)+
\ov{g\left(u(t),\ov u(t)\right)}\cdot\phi(0)=$$
$$=Aw\left(0,u(t),\ov u(t)\right)+\al Gw\left(-\tau_0, u(t),\ov u(t)\right)+
N\left(u(t)\phi(0)+\ov u(t)\ov\phi(0)+\right.$$ $$+\left.w\left(0,
u(t), \ov u(t)\right), u(t)\phi(-\tau_0)+\ov
u(t)\ov\phi(-\tau_0)+w\left(-\tau_0, u(t),\ov
u(t)\right)\right),\eqno(6.33)$$ where $u(t)$ is a solution of the
ordinary differential equation
$$\dot u(t)=\lam_1 u(t)+g\left(u(t),\ov u(t)\right)\eqno(6.34)$$
and
$$g\left(u,\ov u\right)=\ov{\widt\psi}(0)N\left(u\phi(0)+\ov u\ov\phi(0)+w\left(0, u,\ov u\right),
u\phi(-\tau_0)+\ov u\ov\phi(-\tau_0)+\right.$$
$$+\left.w\left(-\tau_0, u,\ov u\right)\right).\eqno(6.35)$$

For $\fii^c\in W_\loc^c(0)\cap V_1$ the solution of (6.20) is
given by
$$V(t)(\th)=u(t)\phi(\th)+\ov u(t)\ov\phi(\th)+w\left(\th, u(t),\ov u(t)\right),\quad\th\in[-\tau_0, 0].\eqno(6.36)$$
Because the equation (6.20) has a Casimir function (conservation
laws) $\|IV\|=m$, for the solution given by (6.36) we have
$IV(t)\cdot I\dot V(t)=0$.

Consider the function $w$ given by
$$w\left(\th, u,\ov u\right)=\f{1}{2}w_{20}(\th)u^2+w_{11}(\th)u\ov u+\f{1}{2}w_{02}(\th)u^2,\eqno(6.37)$$
with $w_{02}(\th)=\ov w_{20}(0)$, $w_{11}(\th)=\ov w_{11}(\th)$,
$\th\in[-\tau_0,0]$. From (6.20) it results that the components of
$N\left(V,\widt V\right)$ for $V=(x^1,x^2,x^3)^T$, $\widt
V=\left(\widt x{}^1,\widt x{}^2,\widt x{}^3\right)^T$ are the
followings:
$$N^1\left(V,\widt V\right)=\f{I_2-I_3}{I_1} x^2x^3+\al m\left[\f{I_2(I_1-I_2)}{I_1}
x^2\widt x{}^2 -\f{I_3(I_3-I_1)}{I_1}x^3\widt x{}^3\right]+$$
$$+\al\left[\f{I_2(I_1-I_2)}{I_1}x^1x^2\widt x{}^2-\f{I_3(I_3-I_1)}{I_1}\widt x{}^1 x^3\widt x{}^3\right],$$
$$N^2\left(V,\widt V\right)=\f{I_3-I_1}{I_2}x^1x^3+\al m\f{I_1(I_2-I_1)}{I_2}
\left(x^1+\widt x{}^1\right)\widt x{}^2+$$
$$+\al\left[\f{I_3(I_2-I_3)}{I_2}\widt x{}^2 x^3\widt x{}^3-\f{I_1(I_1-I_2)}{I_2}x^1
\widt x{}^1\widt x{}^2\right],\eqno(6.38)$$
$$N^3\left(V,\widt V\right)=\f{I_1-I_2}{I_3}x^1x^2+\al m\f{I_1(I_3-I_1)}{I_3}\left(x^1+\widt x{}^1\right)\widt x{}^3+$$
$$+\al\left[\f{I_1(I_3-I_1)}{I_3} x^1\widt x{}^1\widt x{}^3-\f{I_2(I_2-I_3)}{I_3}x^2
\widt x{}^2\widt x{}^3\right].$$

From (6.38) with $V(t)$, $\widt V(t)=V(t)(-\tau_0)$, given by
(6.36) and $\phi(\th)$, $\ov\phi(\th)$ given by (6.24), it results
$$N\left(V(t),\widt V(t)\right)=\f{1}{2} F_{20}+F_{11} u(t)\ov u(t)+\f{1}{2}
F_{02}\ov u(t)^2+\f{1}{2}F_{21}u(t)^2\ov u(t),\eqno(6.39)$$ with
$F_{20}=(F_{20}^1, F_{20}^2, F_{20}^3)^T$, $F_{11}=(F_{11}^1,
F_{11}^2, F_{11}^3)^T$, $F_{02}=(F_{02}^1, F_{02}^2, F_{02}^3)^T$,
$F_{21}=(F_{21}^1, F_{21}^2, F_{21}^3)^T$, where
$$\begin{array}{l}
\vspace{0.2cm}
F_{20}^1=\ds\f{2(I_2-I_3)}{I_1}v_2v_3+\ds\f{\al m}{I_1}[I_2(I_1-I_2)v_2^2-I_3(I_3-I_1)v_3^2]e^{\lam_2\tau_0},\\
\vspace{0.2cm}
F_{20}^2=F_{20}^3=0,\\
\vspace{0.2cm} \begin{array}{lll} \vspace{0.2cm} \!\!\!
F_{11}^1&=&\ds\f{I_2-I_3}{I_1}v_2\left(v_3+\ov
v_3\right)+\ds\f{\al m}{I_1}\left[I_2(I_1-I_2)v_2^2
-\right.\\
&-&I_3(I_3-I_1)v_3\ov v_3\Bigl](e^{\lam_1\tau_0}+e^{\lam_2\tau_0}),\end{array}\\
\vspace{0.2cm}
F_{11}^2=F_{11}^3=0,\\
\vspace{0.2cm}
F_{02}^1=\ds\f{2(I_2-I_3)}{I_1}v_2\ov v_3+\ds\f{\al m}{I_1}\left[I_2(I_1-I_2)v_2^2-I_3(I_3-I_1)\ov v{}_3^2\right] e^{\lam_1\tau_0},\\
\vspace{0.2cm} F_{02}^2=F_{02}^3=0,\\\vspace{0.2cm}
 \begin{array}{lll} \vspace{0.2cm}
\!\!\!F_{21}^1\!\!\!&=&\!\!\ds\f{I_2-I_3}{I_1}\left[v_2(2w_{11}^3(0)+w_{20}^3(0))+2v_3w_{11}^2(0)+\ov v_3 w_{20}^2(0)\right]+\\
\vspace{0.2cm}
&+&\!\!\ds\f{2\al m I_2(I_1\!-\!I_2)}{I_1}\left[2v_2w_{11}^2(0)e^{-\lam_1\tau_0}\!+\ds\f{1}{2} v_2w_{20}^2(e^{\lam_1\tau_0}\!+e^{\lam_2\tau_0})\right]-\\
&-&\!\! \ds\f{2\al m I_3(I_3\!-\!I_1)}{I_1}\left[v_3 w_{11}^3(0)e^{-\lam_1\tau_0}\!+\ds\f{1}{2}\ov v_3 w_{20}^3(0)(e^{\lam_1\tau_0}\!+e^{\lam_2 \tau_0})\right]\!,\\
\end{array}\\
\vspace{0.2cm} F_{21}^2=\ds\f{I_3-I_1}{I_2}\ov v_3 w_{20}^1(0)-\ds\f{2\al m
I_1(I_1-I_2)}{I_2} v_2 w_{20}^1(0)e^{\lam_1\tau_0},\end{array}\eqno(6.40)$$
$$F_{21}^3=\ds\f{I_1-I_2}{I_3} v_2 w_{20}^1(0)-\ds\f{2\al m I_1(I_1-I_3)}{I_3}\ov v_3
w_{20}^1(0)e^{\lam_1\tau_0}.$$

From (6.35) with $\widt\psi(0)$ given by (6.28) and
$N\left(V(t),\widt V(t)\right)$ given by (6.39) we obtain
$$g\left(u(t),\ov u(t)\right)=\f{1}{2} g_{21} u(t)^2\ov u(t),\eqno(6.41)$$
where
$$g_{21}=\ov{\widt w}_2F_{21}^2+\ov{\widt w}_3 F_{21}^3.\eqno(6.42)$$

Taking into account of (6.32) it results that $w_{20}(\th)$,
$w_{11}(\th)$ verify the differential equations
$$\dot w_{20}(\th)=2\lam_1 w_2(\th),\quad \dot w_{11}(\th)=0,\quad \th\in [-\tau_0, 0].\eqno(6.43)$$

From (6.43), (6.33) and $IV(t)\cdot I\dot V(t)=0$ it results
$$w_{20}(\th)=E_1e^{2\lam_1\th},\quad w_{11}(\th)=0,\quad \th\in [-\tau_0, 0],\eqno(6.44)$$
where $E_1$ is the solution of the linear system of equations
$$\left(A+\al e^{-\lam_1\tau_0}G-2\lam_1E\right)E_1=-F_{20}.\eqno(6.45)$$
From (6.45) and (6.40) we deduce
$$w_{20}^1(\th)=\f{1}{2\lam_1}F_{20}^1 e^{2\lam_1\th},\quad w_{20}^2(\th)=w_{20}^3(\th)=0,\quad \th\in [-\tau_0,0],\eqno(6.46)$$
obtaining the following

{\bf Proposition 6.4.} {\it The solution of the equation $(6.6)$
near upon the stationary state
$\Om_1=\left(\ds\f{m}{I_1},0,0\right)$ is
$$\begin{array}{l}
\vspace{0.2cm}
\dot x(t)=\ds\f{m}{I_1}+\Rem(w_{20}^1(0)u^2(t)),\\
\vspace{0.2cm}
\dot y(t)=2 v_2\Rem (u(t)),\\
\dot z(t)=2\Rem (v_3 u(t)),\end{array}\eqno(6.47)$$ where $u(t)$
is a solution of the equation
$$\dot u(t)=\f{1}{2}g_{21}u(t)^2 \ov u(t)\eqno(6.48)$$
and $w_{20}^1(0)=\ds\f{1}{2\lam_1}F_{20}^1$, $v_2, v_3$ are given
by $(6.24)$, $g_{21}$ is given by $(6.42)$.}

Based on the above analysis and calculi, we can see that $g_{21},
w_{20}^1(0), v_2, v_3$ are determined by the parameters and the
delay of (6.6). Thus we can explicitly compute the following
quantities:
$$\begin{array}{l}
\vspace{0.2cm} C_1(0)=\ds\f{1}{2}g_{21},\quad
\mu_2=-\ds\f{\Rem(C_1(0))}{\Rem\left(\f{d\lam}{d\tau}\right)_
{\lam=i\om_0\atop{\tau=\tau_0}}},\\
T_2=-\ds\f{Im(C_1(0))+\mu_2 Im\left(\f{
d\lam}{d\tau}\right)_{\lam=i\om_0\atop{\tau=\tau_0}}}{\om_0},~\beta_2=2\Rem(C_1(0)).\end{array}
\eqno(6.49)$$

In summary, this leads to the following result.

{\bf Proposition 6.5.} {\it In the formulas $(6.49)$ $\mu_2$
determines the direction of the Hopf bifurcation: if $\mu_2>0$
(respectively $\mu_2<0$) then the Hopf bifurcation is
supercritical (respectively subcritical) and the bifurcating
periodic solutions exist for $\tau>\tau_0$ (respectively
$\tau<\tau_0$); $\beta_2$ determines the stability of the
bifurcation periodic solutions: the solutions are orbitally stable
(respectively instable) if $\beta_2<0$ (respectively $\beta_2>0$);
$T_2$ determines the periods of the bifurcating periodic
solutions: the periods increase (respectively decrease) if $T_2>0$
(respectively $T_2<0$).} For $I_1=0.8$, $I_2=0.5$, $I_3=0.4$,
$\al=0.3$, $m=1.5$ and $\om_0$, $\tau_0$ given by the formulas
(6.14) we obtain $\om_0= 3.20631$, $\tau_0=0.88154$,
$\mu_2=0.00958$, $T_2=0.00057$, $\beta_2=-0.00139$. The limit
cycle is supercritical with the period $T_2$.  For $I_1=0.8$,
$I_2=0.5$, $I_3=0.4$, $\al=0.3$, $m=1.8$ and $\om_0$, $\tau_0$
given by the formulas (6.15) we obtain $\om_0=0.68547$,
$\tau_0=0.88154$, $\mu_2=0.00344$, $T_2=0.00050$,
$\beta_2=0.00097$. The limit cycle is supercritical with the
period $T_2$.
\section*{\normalsize\bf 7. Conclusions and comments}

\hspace{0.6cm} In this paper we have given a general method to
construct a dissipative mechanism with delay preserving the
symplectic leaves of the reduced space and dissipating the energy.
The most important case is that of the dual of a Lie algebra when
the dissipative term with delay is shown to have a double bracket
with delay form. This theory applies to a number of interesting
examples from ferromagnetics, ideal fluid flow and plasma dynamics
in which the previous state of the phenomenon is important. In the
future we would like to analyze the systems given by
Landau--Lifschitz equations with delay, the Heavy Top with delay
etc.

\noindent{\footnotesize{\begin{tabular}{lllll}
I.D. Albu & M. Neam\c tu & D. Opri\c s\\
Dept. of Mathematics & Faculty of Economics & Dept. of Appl. Math.\\
West Univ. of Timi\c soara & West Univ. of Timi\c soara & West Univ. of Timi\c soara\\
albud@math.uvt.ro & mihaela.neamtu@fse.uvt.ro & opris@math.uvt.ro
\end{tabular}
}}

\end{document}